\documentclass[10pt,english]{smfart}
\usepackage[T1]{fontenc}
\usepackage{lmodern}
\IfFileExists{french.ldf}{\usepackage[french,english]{babel}}{%
  \usepackage[english]{babel}\babelprovide[import]{french}}
\usepackage{amsmath,amsfonts,amssymb,mathrsfs}
\usepackage{enumerate,url,xspace,smfthm}
\IfFileExists{smfhyperref.sty}{\usepackage[hidelinks]{hyperref}}{}
\newcommand{\F}{\mathbb F}
\newcommand{\A}{\mathbb A}
\newcommand{\Pj}{\mathbb P}

\newcommand{\Ocal}{\mathcal O}
\newcommand{\m}{\mathfrak m}
\newcommand{\p}{\mathfrak p}
\newcommand{\kbar}{\overline{\F}_q}
\DeclareMathOperator{\Sing}{Sing}
\DeclareMathOperator{\Supp}{Supp}
\DeclareMathOperator{\Spec}{Spec}
\DeclareMathOperator{\Proj}{Proj}

\DeclareMathOperator{\ord}{ord}
\DeclareMathOperator{\gr}{gr}
\DeclareMathOperator{\length}{length}

\numberwithin{equation}{section}
\title[Counting multiplicities in projective varieties]{Counting multiplicities in projective varieties over a finite field}
\alttitle{Comptage des multiplicit\'es dans les vari\'et\'es projectives sur un corps fini}
\author{Chunhui Liu}
\address{Institute for Advanced Study in Mathematics\\
Harbin Institute of Technology\\150001 Harbin\\P. R. China}
\email{chunhui.liu@hit.edu.cn}
\date{\today}

\begin{document}
\raggedbottom
\def\smfbyname{}

\begin{abstract}
In this article, we fix a counting function, and consider a problem of counting multiplicities in a pure-dimensional projective variety over a finite field. We give an optimal upper bound of the sum with respect to this counting function, which gives a description of the global complexity of its singular locus. 

\end{abstract}
\begin{altabstract}
Dans cet article, on fixe une fonction de comptage et \'etudie un problème de comptage des multiplicit\'es sur une vari\'et\'e projective de dimension pure d\'efinie sur un corps fini. On \'etablit une majoration optimale de la somme associ\'ee \`a cette fonction de comptage, ce qui permet de d\'ecrire la complexit\'e globale du lieu singulier de cette vari\'et\'e.
\end{altabstract}

\renewcommand{\subjclassname}{2020 Mathematics Subject Classification}
\maketitle
\tableofcontents

\section{Introduction}\label{sec:intro}

\subsection{Rational points and local multiplicities}
The problem of counting rational points over a finite field asks how
the number of solutions of polynomial equations is constrained by
their geometry. For a projective scheme, its dimension determines
the natural power of the cardinality of the field, while its degree
provides a measure of its geometric complexity. A basic estimate in
this direction is
\[
 \#X(\F_q)\leqslant\delta\pi_d(q),
 \qquad \pi_d(q)=1+q+\cdots+q^d,
\]
for a reduced projective scheme $X$ of pure dimension $d$ and degree
$\delta$ over the finite field $\F_q$ of cardinality $q$
\cite[Proposition~2.17, pp.~331--332]{Liu2022}. The polynomial $\pi_d(q)$ is the
number of rational points of $\Pj^d_{\F_q}$. More generally, we use
the notation $\pi_j(q)=1+q+\cdots+q^j$ for every integer $j\geqslant0$.

When singularities are present, the number of rational points alone
does not describe their contribution to the geometry of $X$.
Even the number of singular rational points does not distinguish
points with different multiplicities. To take this information into
account, one may attach a weight to each rational point according to
the multiplicity of its local ring. We are interested in uniform
upper bounds for the resulting sums, expressed in terms of the
degree of $X$, the dimension of its singular locus, and the
cardinality of the field.

\subsubsection{}
Let $X\subset\Pj^n_{\F_q}$ be a nonempty reduced closed subscheme of
pure dimension $d$ and degree $\delta$. For a point $x$ of $X$, we
denote by $\mu_x(X)$ the Hilbert--Samuel multiplicity of the local
ring $\Ocal_{X,x}$ with respect to its maximal ideal $\m_x$.
The definition and the local properties used below are recalled in
Section~\ref{sec:algebra}. A smooth point has multiplicity one.
Thus a counting function which vanishes at one removes the
contribution of the smooth locus. For a positive integer $t$, we
consider
\begin{equation}\label{eq:moment}
 S_t(X;q)=\sum_{x\in X(\F_q)}\mu_x(X)\bigl(\mu_x(X)-1\bigr)^t.
\end{equation}
The factor $\mu_x(X)$ retains the local multiplicity, while each
factor $\mu_x(X)-1$ increases the weight of points with higher
multiplicity. This choice is also adapted to intersection theory,
as the classical plane-curve estimate recalled below illustrates.

We denote by $\Sing(X)$ the complement of the smooth locus and
put $\dim\varnothing=-1$. Since $\F_q$ is perfect, smoothness is
equivalent to regularity at each point of $X$
\cite[Proposition~6.7.4 and Th\'eor\`eme~6.8.6]{EGAIV2}. The sum~\eqref{eq:moment} combines
the sizes of the local multiplicities with the distribution of the
singular rational points. For example, for every integer $m\geqslant2$,
\begin{equation}\label{eq:intro-multiplicity-tail}
 \#\{x\in X(\F_q):\mu_x(X)\geqslant m\}
 \leqslant\frac{S_t(X;q)}{m(m-1)^t},
\end{equation}
because every point in this set contributes at least $m(m-1)^t$
to the sum. A bound for~\eqref{eq:moment} therefore controls the
number of points above every multiplicity threshold.

\subsection{Geometric and arithmetic origins}

For plane curves, sums of multiplicities occur naturally in the
classical applications of B\'ezout's theorem. Let $C\subset\Pj^2_k$
be a reduced curve of degree $\delta$ over an algebraically closed
field $k$. One has
\begin{equation}\label{eq:intro-plane-curve}
 \sum_{x\in C(k)}\mu_x(C)\bigl(\mu_x(C)-1\bigr)
 \leqslant\delta(\delta-1).
\end{equation}
See \cite[Exercise~5-22]{FultonCurves} for the irreducible case
and \cite[\S1.1, p.~314]{Liu2022} for the reduced formulation above.
If $\delta\geqslant2$, a general polar curve has degree $\delta-1$
and has order at least
$\mu_x(C)-1$ at $x$. Its intersection multiplicity with $C$ at
$x$ is therefore at least $\mu_x(C)(\mu_x(C)-1)$, and the total
intersection degree is $\delta(\delta-1)$. This explains both
the local weight and the degree of the right side of
\eqref{eq:intro-plane-curve}.

The same relation between global degree and local singularity
invariants appears in the theory of dual varieties. Over an
algebraically closed field of characteristic zero, Laumon computed
in 1976 the degree of the dual of a reduced hypersurface with
isolated singularities, assuming that its dual is a hypersurface
\cite[Corollaire~4.2.1, p.~60]{Laumon1976}. His formula expresses the
difference between the degree of the dual and the corresponding
degree for a smooth hypersurface in terms of local invariants of
the singularities. For a hypersurface of degree $\delta$ in
$\Pj^n$, the latter degree is $\delta(\delta-1)^{n-1}$.
The isolated-singularity case thus suggests estimates which do
not involve the cardinality of the ground field. When the
singular locus has positive dimension, a sum over all geometric
singular points is no longer the appropriate finite quantity.
Over a finite field, summing over rational points gives such a
quantity, and the size of the field must then enter the estimate.

\subsubsection{}
An arithmetic motivation comes from the determinant method for
rational points of bounded height. In this method, rational points
are grouped according to their reductions at finite places, and
one constructs auxiliary hypersurfaces containing the points in
each group. The local estimates depend on the Hilbert--Samuel
functions at the reduction points. In 2007, Salberger made this dependence on the multiplicities
explicit in his criterion for the existence of an auxiliary
hypersurface
\cite[Theorem~3.2]{Salberger2007}. Chen developed in 2012 a formulation
using the slope method in Arakelov geometry, with the local
Hilbert--Samuel functions entering the auxiliary-hypersurface
criterion \cite[Theorem~3.1]{Chen2012II}.

These local estimates make it useful to control how many reduction
points have large multiplicity. Whenever a reduction is reduced
and pure dimensional, an upper bound for~\eqref{eq:moment}
provides precisely this control through
\eqref{eq:intro-multiplicity-tail}. This arithmetic question and
the geometric problem of measuring the singular locus lead to
the same finite-field counting problem. In higher dimension, one
must account simultaneously for the multiplicities and for the
dimension of the locus on which they occur.

\subsection{Bounds over finite fields}

In our 2022 work \cite{Liu2022}, we studied~\eqref{eq:moment} for a reduced
hypersurface $X\subset\Pj^n_{\F_q}$ whose nonempty singular locus
has dimension $s$. In this case $\dim X=n-1$, and the codimension
of the singular locus in $X$ is $n-s-1$. For this exponent,
\cite[Theorem~5.1, p.~346]{Liu2022} gives
\begin{equation}\label{eq:old}
 S_{n-s-1}(X;q)
 \leqslant\sum_{j=0}^s
      \delta(\delta-1)^{n-s-1+j}\pi_{s-j}(q).
\end{equation}
The proof uses intersection trees to organize successive sections
by derivative hypersurfaces and to bound the multiplicities of
the components which arise. When $s=0$, the right side is
$\delta(\delta-1)^{n-1}$. When $s>0$, the additional terms account
for components of smaller dimension, and each decrease of the
dimension is accompanied by a further factor of $\delta-1$.

Estimate~\eqref{eq:old} holds over every finite field. In particular,
it gives
\[
 S_{n-s-1}(X;q)\ll_n
 \delta^{n-s}\max\{\delta-1,q\}^{s}.
\]
Consequently, when $q\geqslant\delta-1$, one obtains the expected
order $\delta^{n-s}q^s$. For fixed $q$ and $s>0$, however, the
terms with $j>0$ in~\eqref{eq:old} grow faster with $\delta$ than
the first term. Moreover, the construction starts from a single
equation defining $X$, and does not directly give the same
estimate for a scheme of higher codimension.

\subsubsection{}
For a reduced closed subscheme $X\subset\Pj^n_{\F_q}$ of pure
dimension $d$, degree $\delta$, and nonempty singular locus of
dimension $s$, \cite[Conjecture~5.13, p.~354]{Liu2022} asks whether
\[
 S_{d-s}(X;q)\ll_n\delta^{d-s+1}q^s.
\]
Here the implied constant is to be uniform in $q$, $X$, and its
degree. The exponent $d-s$ is the codimension of the singular
locus in $X$. It is also the number of proper hypersurface
sections needed to reach its dimension, which explains its role
in an estimate obtained by intersection theory.

We prove this conjecture by establishing the sharper bound
$\delta(\delta-1)^{d-s}\pi_s(q)$, with no restriction on $q$.
The result applies in arbitrary codimension. More generally, it
gives a bound for every positive exponent $t$ up to $d-s$, with
right side $\delta(\delta-1)^t\pi_{d-t}(q)$. Thus it both extends
the hypersurface estimate and removes its additional terms.
Cones over smooth varieties show that the resulting coefficient
one is attained over every finite field. We now state the theorem
and its consequences precisely.

\subsection{Main results}

We prove a family of estimates in which the allowed exponents are
determined by the codimension of the singular locus. The bound is
uniform in the degree and in the finite field, and its coefficient
does not depend on the ambient dimension.

\begin{theo}\label{thm:main}
Let $X\subset\Pj^n_{\F_q}$ be a nonempty reduced closed subscheme of
pure dimension $d\geqslant1$ and degree $\delta$.
For every integer $t$ such that $1\leqslant t\leqslant d$ and
$\dim\Sing(X)\leqslant d-t$, one has
\begin{equation}\label{eq:family}
 S_t(X;q)\leqslant\delta(\delta-1)^t\pi_{d-t}(q).
\end{equation}
\end{theo}

The two factors on the right of~\eqref{eq:family} arise separately
in the proof. Intersecting $X$ properly with $t$ hypersurfaces of
degree $\delta-1$, when $\delta>1$, gives a cycle of degree
$\delta(\delta-1)^t$. Its dimension is $d-t$, which accounts for the
factor $\pi_{d-t}(q)$ in the point count. If $\delta=1$, $X$ becomes
a linear space over the algebraic closure by
Proposition~\ref{prop:finite-field-transfer} and Lemma~\ref{lem:degree-one}, and both sides of~\eqref{eq:family} vanish.

Pure dimensionality is needed here. With the usual convention that
degree is determined by the leading coefficient of the Hilbert
polynomial, components of smaller dimension can contribute to
\eqref{eq:moment} without contributing to the degree of $X$.
An example is given in Remark~\ref{rem:pure-needed}.

\subsubsection{}
The largest admissible exponent gives the estimate predicted by the
conjecture. The exponent one gives a bound valid for every reduced
scheme of pure positive dimension.

\begin{coro}\label{cor:fixed}
Suppose that $X$ satisfies the hypotheses of Theorem~\ref{thm:main}
on dimension, degree and reducedness. If $\Sing(X)$ is nonempty
of dimension $s$, put $r=d-s$. Then
\begin{equation}\label{eq:main}
 \sum_{x\in X(\F_q)}\mu_x(X)\bigl(\mu_x(X)-1\bigr)^r
 \leqslant\delta(\delta-1)^r\pi_s(q).
\end{equation}
\end{coro}

\begin{proof}
At the generic point of each irreducible component of $X$ the local
ring is a field, since $X$ is reduced. Hence all these generic points
are regular, and therefore smooth over $\F_q$. The smooth locus is open
\cite[Corollaire~6.8.7]{EGAIV2}, so its complement has dimension at
most $d-1$. Thus $1\leqslant r\leqslant d$, and taking $t=r$ in
Theorem~\ref{thm:main} gives~\eqref{eq:main}.
\end{proof}

When $X$ has only isolated singularities, the right side
of~\eqref{eq:main} is independent of $q$. At the other extreme, if
the singular locus has codimension one in $X$, the exponent in
\eqref{eq:main} is one. An estimate for this exponent remains
available when the dimension of the singular locus is unknown.

\begin{coro}\label{cor:reduced}
For every nonempty reduced $X\subset\Pj^n_{\F_q}$ of pure dimension
$d\geqslant1$ and degree $\delta$, one has
\begin{equation}\label{eq:quadratic}
 \sum_{x\in X(\F_q)}\mu_x(X)\bigl(\mu_x(X)-1\bigr)
 \leqslant\delta(\delta-1)\pi_{d-1}(q).
\end{equation}
In particular,
\begin{equation}\label{eq:quadratic-simple}
 S_1(X;q)\leqslant(2-2^{-(d-1)})\delta(\delta-1)q^{d-1}
            \leqslant2\delta^2q^{d-1}.
\end{equation}
If $X$ is reduced of dimension zero, then $S_t(X;q)=0$ for every
$t\geqslant1$.
\end{coro}

\begin{proof}
The preceding proof gives $\dim\Sing(X)\leqslant d-1$ whenever
$d\geqslant1$. We take $t=1$ in Theorem~\ref{thm:main}. To obtain
\eqref{eq:quadratic-simple}, we use the inequality
\begin{equation}\label{eq:intro-pi-estimate}
 \frac{\pi_j(q)}{q^j}=\sum_{i=0}^j q^{-i}
 \leqslant\sum_{i=0}^j2^{-i}=2-2^{-j}
 \qquad(j\geqslant0).
\end{equation}
If $d=0$, the local ring at a rational point is a reduced Artinian
local ring. Its maximal ideal is nilpotent and hence zero. Its
multiplicity is therefore one.
\end{proof}

\subsubsection{}
To compare with the conjectured order of growth, it remains to bound
the ratio $\pi_s(q)/q^s$. The preceding geometric series yields an
absolute constant, even when the degree exceeds the field size.

\begin{coro}\label{cor:conjecture}
Under the hypotheses of Corollary~\ref{cor:fixed}, one has
\begin{equation}\label{eq:big-o}
 S_{d-s}(X;q)
 \leqslant (2-2^{-s})\delta^{d-s+1}q^s
 \leqslant 2\delta^{d-s+1}q^s.
\end{equation}
In particular, Conjecture~5.13 of \cite{Liu2022} holds with an
absolute constant.
\end{coro}

\begin{proof}
Apply~\eqref{eq:intro-pi-estimate} with $j=s$ in~\eqref{eq:main}
and use $(\delta-1)^{d-s}\leqslant\delta^{d-s}$.
\end{proof}

\begin{coro}\label{cor:hypersurface}
Let $X\subset\Pj^n_{\F_q}$ be a reduced hypersurface of degree
$\delta$ with nonempty singular locus of dimension $s$. Then
\begin{equation}\label{eq:hypersurface}
 \sum_{x\in X(\F_q)}\mu_x(X)\bigl(\mu_x(X)-1\bigr)^{n-s-1}
 \leqslant\delta(\delta-1)^{n-s-1}\pi_s(q).
\end{equation}
For every reduced hypersurface in $\Pj^n_{\F_q}$ with $n\geqslant2$,
one also has
\[
 \sum_{x\in X(\F_q)}\mu_x(X)\bigl(\mu_x(X)-1\bigr)
 \leqslant\delta(\delta-1)\pi_{n-2}(q).
\]
\end{coro}

\begin{proof}
The irreducible factors of a homogeneous equation for $X$ define
its irreducible components, each of dimension $n-1$. The two
assertions are Corollaries~\ref{cor:fixed} and~\ref{cor:reduced}
with $d=n-1$.
\end{proof}

The right side of~\eqref{eq:hypersurface} is the term $j=0$
in~\eqref{eq:old}. Thus the estimates agree when $s=0$, and
\eqref{eq:hypersurface} removes all the terms with $j>0$ when $s>0$.
The remaining term is attained by cones. Indeed, for a cone of
degree $\delta>1$ over a smooth base with vertex $\Pj^s_{\F_q}$,
each rational point of the vertex has multiplicity $\delta$, and the
other points are smooth. This also explains the factor $\pi_s(q)$
in~\eqref{eq:main}. Section~\ref{sec:sharpness} proves this assertion,
extending \cite[Example~5.11, p.~353]{Liu2022}, and gives an equality
example in codimension two. Corollary~\ref{cor:optimal-n} determines
the least constant for each fixed ambient dimension in the normalization
of~\eqref{eq:big-o}.

\subsection{Method of proof}

The plane-curve estimate suggests intersecting with polar
hypersurfaces and comparing the local contributions with the total
degree. In higher dimension the intersections have positive dimension,
so a point count is needed as well. In arbitrary codimension one must
first replace the single defining equation of a hypersurface by a
family of equations with the same local multiplicity properties.

\subsubsection{}
Suppose first that $X$ is a hypersurface with equation $F$.
At a geometric point $x$, the order of a local equation is
$\mu_x(X)$. A first derivative has order at least $\mu_x(X)-1$,
and at a smooth point some first derivative is nonzero. Thus the
partial derivatives give forms of degree $\delta-1$ which vanish
to the required order at singular points and whose common zero set
on $X$ is its singular locus.

For a reduced scheme of higher codimension we use finite linear
projections to $\Pj^{d+1}$ over $\overline{\F}_q$. We take the
image with its cycle multiplicities: the image of each irreducible
component is counted with the degree of the corresponding extension
of function fields, and components with the same image are counted
together. Pulling back an equation of this image cycle gives a form
of degree $\delta$ which vanishes on $X$. Its local order at every
geometric closed point is at least the multiplicity of $X$ there.
At a prescribed smooth point the projection can be chosen so that
the pulled-back equation has a nonzero linear term. These are the
properties of Chow equations needed here, in the projective setting
of Koll\'ar's construction
\cite[Definition~4.1 and Example~4.6]{Kollar1999}.
Section~\ref{sec:projection} identifies these equations on the affine
cone and transfers Koll\'ar's local statements to projective space.
The degree calculation includes inseparable projections and reducible
images.

\subsubsection{}
Taking first derivatives of these equations gives a system of forms
of degree $\delta-1$ whose base locus on
$X_{\overline{\F}_q}$ is contained in the singular locus. We choose
$t$ successive members of this system. Before the $i$th section,
every component of the current cycle has dimension $d-i+1>d-t$.
The hypothesis $\dim\Sing(X)\leqslant d-t$ therefore ensures that
each such component contains a smooth point of $X_{\overline{\F}_q}$.
Over the algebraic closure one can choose a member which contains
none of the finitely many components. The section is consequently
proper.

At a point of multiplicity greater than one, the local section
inequality supplies a factor of $\mu_x(X_{\overline{\F}_q})-1$
at each step. The resulting effective cycle $C$ has pure dimension
$d-t$, and satisfies
\[
 \deg C=\delta(\delta-1)^t,\qquad
 \mu_x(C)\geqslant
 \mu_x(X_{\overline{\F}_q})
 \bigl(\mu_x(X_{\overline{\F}_q})-1\bigr)^t.
\]
The inequality is also valid at smooth points, since its right side
is zero. Sections~\ref{sec:algebra} and~\ref{sec:main-proof}
apply the local inequality to cycles and carry out this successive
construction.

The remaining step is to count with the multiplicities of the cycle
still attached. Section~\ref{sec:count} proves
\[
 \sum_{x\in\Pj^n(\F_q)}\mu_x(C)
 \leqslant (\deg C)\pi_{d-t}(q).
\]
Here the points belong to the fixed ambient space over $\F_q$,
although the cycle is allowed to be defined only over its algebraic
closure. This is the step that avoids further decomposing the
singular locus into pieces of different dimensions, as in the
intersection-tree estimate~\eqref{eq:old}. The multiplicity at a
rational point of $X$ is unchanged by extension of the ground field,
so the two displayed inequalities prove Theorem~\ref{thm:main}.

This distinction between the field of definition of the auxiliary
cycle and the field of the points being counted also removes the
restriction on the size of $q$. Neither the centers of projection
nor the hypersurfaces used for the proper sections have to be chosen
over $\F_q$. Finally, Section~\ref{sec:sharpness} computes the
multiplicities of cones and constructs equality cases over every
finite field, including examples which are not local complete
intersections.

\section{Local algebra, degrees, and hypersurface sections}\label{sec:algebra}

We first recall the notions of multiplicity and degree which enter the
counting problem. They measure the leading terms of two Hilbert
polynomials, one attached to a local ring and the other to a homogeneous
coordinate ring. Both invariants can be written as sums over
irreducible components, with their generic lengths as coefficients. This allows us to keep track of the multiplicities which
appear under successive intersections.

For a proper hypersurface section, the degree is multiplied by the
degree of the hypersurface. At a point of the section, the new local
multiplicity is at least the old multiplicity times the vanishing order
of the equation. We recall the global
degree formula and the standard local multiplicity inequality. Successive sections need not be Cohen--Macaulay, so their
fundamental cycles will be used throughout.

\subsection{Local multiplicities}

We recall the Hilbert--Samuel multiplicity and the associativity
formula before passing to the geometric setting. The ideal is
specified in the notation whenever it differs from the maximal ideal.

\subsubsection{}
For a module $M$ of finite length over a ring $A$, we denote its length
by $\length_A(M)$. It is the number of simple quotients in a composition
series and is additive in short exact sequences
\cite[Chap.~6]{AtiyahMacdonald}. Over a local ring $(A,\m)$ all
simple modules are isomorphic to the residue field $A/\m$. In
particular, if $A$ is a local algebra over a field $k$ with residue
field $k$, the length of a finite-length module is its dimension as a
$k$-vector space. We use this observation when extending the ground field.

For a finite $A$-module $M$, its support is
$\Supp(M)=\{\p\in\operatorname{Spec}(A):M_\p\ne0\}$, and its
dimension is the dimension of this support. If $(A,\m)$ is a
Noetherian local ring, $I$ is an $\m$-primary ideal,
and $M$ is a nonzero finite $A$-module of dimension $t$, then
$M/I^{N+1}M$ has finite length. Here $I$ is $\m$-primary precisely
when $\sqrt I=\m$. The function $\length_A(M/I^{N+1}M)$ agrees,
for sufficiently large integers $N$, with a polynomial of degree $t$.
Its leading term is
\begin{equation}\label{eq:HS}
 \length_A(M/I^{N+1}M)
 =\frac{e_I(M)}{t!}N^t+O(N^{t-1})\qquad(t\geqslant1).
\end{equation}
For $t=0$, the module $M$ has finite length and $I^{N+1}M=0$ for
large $N$, so $e_I(M)=\length_A(M)$. For the existence and degree
of the polynomial, see \cite[\S13, Theorems 13.2 and 13.4]{Matsumura}.
The coefficient $e_I(M)$ is a positive integer
\cite[\S14, pp.~107--108]{Matsumura}, called the Hilbert--Samuel multiplicity.
We abbreviate $e_\m(A)$ to $e(A)$.
The relation with the graded Hilbert function is given by
\[
 \length_A(M/I^{N+1}M)
 =\sum_{j=0}^N\length_A(I^jM/I^{j+1}M).
\]
Thus, for $t\geqslant1$, the eventual polynomial of the summand
has leading term $e_I(M)j^{t-1}/(t-1)!$. Both conventions give the
same multiplicity.

Let $0\to M'\to M\to M''\to0$ be a short exact sequence of finite
modules of dimension at most $t$. The coefficients of $N^t$ in their
Hilbert--Samuel polynomials are additive, with a zero contribution in
dimension less than $t$ \cite[Theorem 14.6]{Matsumura}. Applied over
$A/\operatorname{Ann}(M)$, the associativity formula
\cite[Theorem 14.7]{Matsumura} gives
\begin{equation}\label{eq:associativity}
 e_I(M)=\sum_{\substack{\p\in\Supp(M)\\\dim A/\p=t}}
 \length_{A_\p}(M_\p)\,e_I(A/\p),
 \qquad t=\dim M,
\end{equation}
where $I$ acts on $A/\p$ by its
image. The primes in the sum are minimal in $\Supp(M)$, so the
localized modules occurring in the coefficients have finite length.

\subsubsection{}
In the computations with homogeneous equations we use unique
factorization in polynomial rings over a field
\cite[Exercises 20.1--20.2]{Matsumura}. We use Nakayama's lemma
in the form of \cite[Proposition 2.6]{AtiyahMacdonald} and the
Nullstellensatz in the form of \cite[Chap.~I, Theorem 1.3A]{Hartshorne}.

For a prime $\p$ in a finitely generated domain $R$ over a field,
the dimension formula gives
$\operatorname{ht}(\p)+\dim(R/\p)=\dim R$
\cite[Chap.~I, Theorem 1.8A]{Hartshorne}.
In particular, a closed point on an integral $t$-dimensional variety
has local dimension $t$, and a nonzero nonunit equation on such a
variety has irreducible zero-set components of dimension $t-1$.
The last assertion also uses the principal ideal theorem
\cite[Theorem~13.5]{Matsumura}.

\subsubsection{}
A local ring $(A,\m)$ with residue field $k$ is regular if
$\dim_k(\m/\m^2)=\dim A$. The vector space $\m/\m^2$ is the
cotangent space at its closed point. If $A$ is regular of dimension
$t$, then $\gr_\m A\simeq k[Z_1,\ldots,Z_t]$
\cite[Theorem 14.4]{Matsumura}. Here
$\gr_\m A=\bigoplus_{j\geqslant0}\m^j/\m^{j+1}$.
In particular $e(A)=1$ \cite[\S14, p.~108]{Matsumura}.
This gives the assertion about smooth points used in the introduction.

\subsection{Cycles and degrees}\label{subsec:cycles}

A hypersurface section of a reduced variety need not be reduced.
We therefore retain each irreducible component with its generic
multiplicity. An effective cycle records exactly these data. Its
multiplicity at a closed point also includes the local multiplicity
of every component passing through that point.

\subsubsection{}
Until Section~\ref{sec:main-proof}, $k$ denotes an algebraically closed
field unless explicitly stated otherwise. An effective pure
$t$-dimensional cycle is a finite sum
$C=\sum_i a_i[V_i]$, where the $V_i$ are distinct integral closed
$t$-dimensional subvarieties and $a_i$ are positive integers. Its
support is the union of the $V_i$. At a closed point $x$, set
\begin{equation}\label{eq:cycle-mu}
 \mu_x(C)=\sum_{i:x\in V_i}a_i\,e(\Ocal_{V_i,x}),
 \qquad \deg C=\sum_i a_i\deg V_i.
\end{equation}
For a component $V$ through $x$, we also write $\mu_x(V)=e(\Ocal_{V,x})$.
The value of $\mu_x(C)$ is zero off the support. The zero element of the group of $t$-cycles has degree and
multiplicities zero.
For an affine variety, degree means the degree of its projective
closure in the given embedding. In this paper, pure dimensionality
means that all irreducible components of the support have the same
dimension. A nonreduced intersection scheme of pure dimension may
still have embedded associated points.

For a projective scheme $V$ of dimension $t$ with standard graded
homogeneous coordinate ring $R$, its degree is defined by
\begin{equation}\label{eq:degree}
 \dim_k R_N=\frac{\deg V}{t!}N^t+O(N^{t-1})\quad(t\geqslant1).
\end{equation}
In dimension zero the eventual constant is $\deg V$. The existence
and degree of this Hilbert polynomial follow from
\cite[Chap.~I, Theorem 7.5]{Hartshorne}.

\subsubsection{}
For a pure-dimensional scheme $V$, its fundamental cycle $[V]$ has
coefficient $\length(\Ocal_{V,\eta})$ at each generic point $\eta$.
The component formulas proved in \cite[Proposition~2.3, p.~321]{Liu2022}
give
\begin{equation}\label{eq:fundamental}
 \deg[V]=\deg V,\qquad \mu_x([V])=e(\Ocal_{V,x})
\end{equation}
for a projective $V$ and a closed point $x\in V$. The local identity
also applies to an affine scheme by taking its projective closure.
Thus the coefficients in the fundamental cycle are the same ones
which enter both the global degree and the local multiplicity.

Formula \eqref{eq:fundamental} remains valid in the presence of
embedded associated points. In the associativity formula at a closed
point, their primes have quotient dimension strictly smaller than
$\dim V$ and therefore do not contribute. For reduced pure-dimensional $V$, every generic length
is one, since a reduced zero-dimensional local ring is a field.
In dimension zero, every component over $k$ is a point of degree one,
so the degree of a scheme is the sum of the lengths of its local rings.
For an effective zero-cycle, its degree is consequently the sum of
its multiplicities over all points of its support.

\subsection{Degree of a hypersurface section}
Let $V$ be an integral subvariety of positive dimension, and let
$H$ be a hypersurface which does not contain $V$. We write
$[V]\cdot H$ for the fundamental cycle of the scheme-theoretic
intersection. More explicitly, if $Z$ runs over the irreducible
components of $V\cap H$ and $\eta_Z$ is its generic point, then
\begin{equation}\label{eq:section-cycle}
 [V]\cdot H
 =\sum_Z\length_{\Ocal_{V,\eta_Z}}
       \bigl(\Ocal_{V,\eta_Z}/(g)\bigr)[Z],
\end{equation}
where $g$ is a local equation of $H$ at $\eta_Z$. The quotient has
finite length because $\eta_Z$ is a minimal prime of the section.
Changing the local equation multiplies $g$ by a unit and leaves the
coefficient unchanged. If $H$ contains no component of
$C=\sum_i a_i[V_i]$, set
$C\cdot H=\sum_i a_i([V_i]\cdot H)$. We say in this case that
the section of $C$ by $H$ is proper.

\begin{lemm}[Degree of a proper hypersurface section]\label{lem:bezout}
Let $C$ be an effective pure $t$-cycle in $\Pj^N_k$, where $t\geqslant1$,
and let $H=V(g)$ be a hypersurface of degree $b\geqslant1$ containing no
component of $C$. Then the hypersurface intersection cycle $C\cdot H$
is effective and pure of dimension $t-1$, and
\begin{equation}\label{eq:bezout}
 \deg(C\cdot H)=b\deg C.
\end{equation}
For a proper affine hyperplane section the corresponding degree is
at most $\deg C$.
\end{lemm}

\begin{proof}
For an integral component $V$ of $C$, the projective hypersurface
form of B\'ezout's theorem \cite[Chap.~I, Theorem~7.7]{Hartshorne}
gives $\deg([V]\cdot H)=b\deg V$. The intersection coefficients
in this theorem are the local lengths in \eqref{eq:section-cycle}.
They are positive, and every component has codimension one in $V$
by the principal ideal theorem. Thus the cycle in the theorem is
the effective pure-dimensional cycle used here. Summing the degree
identity with the coefficients of $C$ gives \eqref{eq:bezout}.

For an affine hyperplane, intersect the projective closure of $V$
with the homogenized hyperplane. This projective intersection cycle
has degree $\deg V$. Its restriction to the affine chart discards
the components at infinity. For each remaining component the chart
contains its generic point, so restriction preserves its coefficient.
Taking the projective closure of its affine part recovers the original
component and its degree. This proves the affine inequality.
\end{proof}

\subsection{The local section inequality}

The local counterpart of the degree formula is the standard
multiplicity inequality recalled below. It applies to a parameter
whose order is bounded below, and does not require the local ring
to be Cohen--Macaulay.

\begin{lemm}\label{lem:local-cut}
Let $(A,\m)$ be a Noetherian local ring of dimension $t\geqslant1$.
Suppose that $f\in\m^a$, where $a\geqslant1$, and that
$\dim A/(f)=t-1$. Then
\begin{equation}\label{eq:local-cut}
 e(A/(f))\geqslant a\,e(A).
\end{equation}
No Cohen--Macaulay hypothesis is required.
\end{lemm}

\begin{proof}
Choose lifts in $A$ of a system of parameters of $A/(f)$.
Together with $f$ they form a system of parameters of $A$.
Apply \cite[Theorem~14.9, p.~109]{Matsumura} with the ideal
$\m$, the module $A$, and the first parameter $f\in\m^a$.
\end{proof}

\begin{coro}\label{cor:cycle-cut}
Let $C$ be an effective pure positive-dimensional cycle in a smooth
ambient affine or projective space, and let $H=V(g)$ contain no
component of $C$. At a closed point $x\in\Supp(C)\cap H$, suppose that a local
equation of $H$ belongs to $\m_{\mathrm{ambient},x}^{a}$ for an
integer $a\geqslant1$. Then
\begin{equation}\label{eq:cycle-cut}
 \mu_x(C\cdot H)\geqslant a\,\mu_x(C).
\end{equation}
\end{coro}

\begin{proof}
For an integral component $V$ through $x$, the image of the equation in
$A=\Ocal_{V,x}$ lies in $\m_A^a$, since the restriction from the
ambient local ring carries the maximal ideal onto $\m_A$. Its image
is nonzero in the domain $A$ by properness, and it is a nonunit
since $x\in H$. Every component of its zero set through $x$ has
dimension $\dim V-1$. This verifies the dimension hypothesis in
Lemma~\ref{lem:local-cut}.

The minimal primes of $A/(g)$ correspond to those components.
Thus \eqref{eq:associativity} and the definition of the hypersurface
intersection cycle give $e(A/(g))=\mu_x([V]\cdot H)$, even if
$A/(g)$ has embedded associated primes. If two components of $C$
give the same component in the section, its generic
coefficients add in the cycle intersection. Multiplying by the
coefficient of $V$ in $C$ and adding over the components through $x$
proves the assertion.
\end{proof}

The exact multiplicity of a projective hypersurface can be read
from the order of its equation. We recall the form of
\cite[Proposition~5.2, p.~346]{Liu2022} used for the target of a
projection.

\begin{lemm}\label{lem:order}
Let $x\in\Pj^N(k)$, and let $g$ be a nonzero local equation of a
hypersurface through $x$. If
$g\in\m_x^a\setminus\m_x^{a+1}$ in the ambient local ring, then
$\mu_x(V(g))=a$.
\end{lemm}

The local equation is obtained by dehomogenizing a defining form
on a chart containing $x$. Changing the chart multiplies it by a
unit, which changes neither its order nor the local hypersurface.
The statement also applies to hypersurfaces in an affine chart,
by homogenizing their equations.

\section{Counting rational points with cycle multiplicities}\label{sec:count}

We now pass from the local inequalities to a finite sum. The
hypersurfaces used later will be chosen over an algebraic closure,
where a nonempty open set of possible equations always has a
$k$-point. Their intersection cycles need not be defined over
$\F_q$. It is therefore necessary to prove the counting estimate
for cycles over $k=\kbar$, while still summing only over points of
the original ambient space defined over $\F_q$. Every $x\in\Pj^N(\F_q)$ determines a unique
$k$-point $\bar x$ by extension of scalars. For a cycle $C$ over
$k$, we abbreviate $\mu_{\bar x}(C)$ to $\mu_x(C)$ when summing
over these ambient rational points. We use the same convention in
affine space. If a cycle is defined only over $k$, the notation $x\in C\cap\A^N(\F_q)$, or its projective analogue, means
that the point obtained by extension of scalars lies on the support
of $C$.

For a pure-dimensional closed subscheme $V\subset\Pj^N_k$,
the unweighted count of these points is bounded by
$(\deg V)\pi_{\dim V}(q)$
\cite[Definition 2.16 and Proposition 2.17]{Liu2022}. We retain the
multiplicity of each point by taking the fundamental cycle after
each hyperplane section.

\subsection{Affine slices}

A nonconstant coordinate on an integral variety gives $q$ disjoint
slices containing all its ambient $\F_q$-points. Each slice has
smaller dimension and degree at most that of the variety. The local
section inequality allows us to perform this induction with point
multiplicities included.

\begin{prop}\label{prop:affine-count}
For every effective pure $t$-cycle $C$ in $\A^N_k$,
\begin{equation}\label{eq:affine-count}
 \sum_{x\in\A^N(\F_q)}\mu_x(C)\leqslant(\deg C)q^t.
\end{equation}
\end{prop}

\begin{proof}
Both sides are additive in cycles, so we can reduce to an integral
variety $V$ and proceed by induction on $t$. For $t=0$, $V$ is a
single $k$-point, its degree and multiplicity are one, and it
contributes either one or zero according to whether it is rational.

Suppose that $t>0$. At least one of the affine coordinates, say $T_i$,
is nonconstant on $V$. Otherwise the coordinate ring of $V$, generated
by those coordinates, would be $k$ and $V$ would have dimension zero.
For each $c\in\F_q$, put $H_c=(T_i=c)$. As $T_i$ is nonconstant,
no $H_c$ contains $V$. Its section cycle $C_c=[V]\cdot H_c$, viewed
in $H_c\simeq\A^{N-1}_k$, is either zero or pure of dimension $t-1$.
Lemma~\ref{lem:bezout} gives $\deg C_c\leqslant\deg V$. The degree
computed inside $H_c$ has the same value as in $\A^N_k$, because
the inclusion of the projective closure of $H_c$ is linear.

At every point of $V\cap H_c$, the equation $T_i-c$ belongs to the
ambient maximal ideal. With order parameter one,
Corollary~\ref{cor:cycle-cut} gives $\mu_x(C_c)\geqslant\mu_x(V)$.
Every ambient rational point has a unique value $c\in\F_q$ for its
coordinate $T_i$. Deleting that coordinate identifies the rational
points of $H_c$ with $\A^{N-1}(\F_q)$. Applying the induction
hypothesis in each hyperplane gives
\[
 \begin{split}
 \sum_{x\in V\cap\A^N(\F_q)}\mu_x(V)
 &\leqslant\sum_{c\in\F_q}
       \sum_{x\in H_c(\F_q)}\mu_x(C_c)\\
 &\leqslant q^{t-1}\sum_{c\in\F_q}\deg C_c
 \leqslant (\deg V)q^t .
 \end{split}
\]
An empty slice has degree zero and contributes zero to both sums.
At the inductive step only the ambient hyperplanes need to be
defined over $\F_q$.
\end{proof}

\subsection{The projective estimate}

For a projective variety, we use one affine chart and its hyperplane
at infinity. The affine contribution gives the term $q^t$, and the
same argument on the hyperplane gives the remaining terms in
$\pi_t(q)=q^t+\pi_{t-1}(q)$.

\begin{prop}\label{prop:projective-count}
For every effective pure $t$-cycle $C$ in $\Pj^N_k$,
\begin{equation}\label{eq:projective-count}
 \sum_{x\in\Pj^N(\F_q)}\mu_x(C)\leqslant(\deg C)\pi_t(q).
\end{equation}
\end{prop}

\begin{proof}
Again reduce to an integral $V$ by additivity. The case $t=0$ is the
same as in Proposition~\ref{prop:affine-count}. For $t>0$, choose a
coordinate hyperplane $H=(T_i=0)$ that does not contain $V$. Such a
hyperplane exists because the intersection of all coordinate
hyperplanes in projective space is empty. The open part
$V\setminus H\subset\Pj^N\setminus H\simeq\A^N$ has projective
closure $V$ and hence degree $\deg V$, since $V\setminus H$ is
a nonempty open subset of the integral variety $V$. Open restriction
preserves local rings, so
each point in this affine part has the same multiplicity as in $V$.
The coordinates $T_j/T_i$, $j\ne i$, identify the chart with
$\A^N$ over $\F_q$. Proposition~\ref{prop:affine-count} thus bounds
the contribution from this chart by $(\deg V)q^t$.

The cycle $D=[V]\cdot H$ is pure of dimension $t-1$, has degree
$\deg V$, and satisfies $\mu_x(D)\geqslant\mu_x(V)$ on $V\cap H$.
The first two assertions follow from Lemma~\ref{lem:bezout}, and
the inequality follows from Corollary~\ref{cor:cycle-cut} for a
linear equation. The induction hypothesis in $H\simeq\Pj^{N-1}$
gives
\[
 \sum_{x\in V\cap H(\F_q)}\mu_x(V)
 \leqslant\sum_{x\in H(\F_q)}\mu_x(D)
 \leqslant (\deg V)\pi_{t-1}(q).
\]
The rational points in the chart and on $H$ form a partition.
Adding the two bounds and using $q^t+\pi_{t-1}(q)=\pi_t(q)$
proves the claim.
\end{proof}

\begin{rema}\label{rem:no-descent}
When $C$ is reducible, we first write $C=\sum_i a_i[V_i]$ and
apply the estimate to each $V_i$. Its proof uses a coordinate
hyperplane chosen for that component. The final sum
$\sum_i a_i(\deg V_i)\pi_t(q)$ therefore involves the given field
size $q$ even when the equations of the components have coefficients
in a finite extension of $\F_q$.
\end{rema}

\section{Finite projections and equations of degree $\delta$}
\label{sec:projection}

Let $k$ be an algebraically closed field and let
$X\subset\Pj^n_k$ be a nonempty reduced closed subscheme of pure
dimension $d<n$ and degree $\delta$. For a hypersurface, the order
of a defining equation at a point is its multiplicity. To use
hypersurface sections in higher codimension, we need equations of
controlled degree that retain this local information. We use the
Chow equations of \cite[Definition~4.1 and Example~4.6]{Kollar1999}.
We describe their projective form, compute their degrees, and explain
how the cited affine statements give the local properties needed below.

\subsection{Finite linear projections}

Independent linear forms $\ell_0,\ldots,\ell_m$ define a projection
\[
 \rho:\Pj^n_k\dashrightarrow\Pj^m_k,
 \qquad [T]\longmapsto[\ell_0(T):\cdots:\ell_m(T)].
\]
Its center is $\Lambda=V(\ell_0,\ldots,\ell_m)$, a linear space of
dimension $n-m-1$ when $m<n$. If $m=n$, the center is empty and
$\rho$ is a projective linear isomorphism. When
$\Lambda\cap X=\varnothing$, the restriction to $X$ is a morphism.

\begin{lemm}\label{lem:finite-projection}
Let $d\leqslant m\leqslant n$. There are independent linear forms
$\ell_0,\ldots,\ell_m$ whose common projective zero set is disjoint
from $X$. Every such choice defines a finite morphism
$\rho:X\to\Pj^m_k$.
\end{lemm}

\begin{proof}
For $m<n$, put $l=n-m-1$. The incidence variety of pairs
$(x,\Lambda)$ with $x\in X\cap\Lambda$ has dimension at most
$d+l(n-l)$, whereas the Grassmannian of $l$-planes has dimension
$(l+1)(n-l)$. Since $d+l<n$, the planes meeting $X$ form a proper
closed subset of this Grassmannian. Its complement is nonempty over
$k$. The case $m=n$ is immediate.

The restriction of $\rho$ to each integral component is proper.
It is also affine, since the inverse image of the target chart
$U_i\ne0$ is the standard affine open $D_+(\ell_i)$ on that
component. It is therefore finite by
\cite[Chapter~II, Exercise~4.6]{Hartshorne}. The same holds for
$X$: on each such chart, its reduced coordinate ring injects into
the product of the coordinate rings of its finitely many components,
which are finite modules over the target chart ring.
\end{proof}

The degree of a finite projection on a component is its
function-field degree, including any inseparable degree. The
projection formula gives the following relation with projective
degree.

\begin{lemm}\label{lem:degree-projection}
For a finite linear projection as above, let $V$ be an integral
component of $X$, let $Y=\rho(V)$ with its reduced structure, and
put $a=[k(V):k(Y)]$. Then
\begin{equation}\label{eq:degree-projection}
 \deg V=a\deg Y.
\end{equation}
\end{lemm}

\begin{proof}
The pushforward of $[V]$ is $a[Y]$, and
$\rho^*\Ocal_Y(1)=\Ocal_V(1)$. Applying the projection formula
\cite[Proposition~2.5(c) and Example~2.5.2(d)]{Fulton1998} to $d$ hyperplane classes
and taking degrees gives \eqref{eq:degree-projection}.
\end{proof}

We also record the degree-one case, which will allow us to exclude
$\delta=1$ in the construction of the derivative system.

\begin{lemm}\label{lem:degree-one}
A reduced pure-dimensional projective scheme of degree one over an
algebraically closed field is a linear subspace.
\end{lemm}

This is \cite[Chapter~I, Exercise~7.6]{Hartshorne}, applied to its
underlying projective algebraic set and its reduced scheme structure.

\subsection{The image cycle and its equation}

An \emph{admissible projection} is a linear projection to
$\Pj^{d+1}_k$ with center disjoint from $X$. When $n=d+1$, we take
the identity projection. Such projections exist by
Lemma~\ref{lem:finite-projection}. Write $X=\bigcup_i X_i$ for
its integral components. Its fundamental cycle is $[X]=\sum_i[X_i]$,
so the image cycle is
\begin{equation}\label{eq:push-cycle}
 D_\rho=\rho_*[X]=\sum_i a_i[Y_i],\qquad
 Y_i=\rho(X_i),\quad a_i=[k(X_i):k(Y_i)].
\end{equation}
Each $Y_i$ is an integral hypersurface in $\Pj^{d+1}$.
If several components have the same image, their coefficients are
added. These coefficients must be retained: they account for both
multiple components above the same image and the degree of the map
on each component.

Choose an irreducible homogeneous equation $g_i$ for $Y_i$ and put
$G_\rho=\prod_i g_i^{a_i}$. If some $Y_i$ coincide, their
irreducible equations differ only by a scalar, so their exponents
add in this product. Thus $G_\rho$ defines the hypersurface cycle
$D_\rho$. Lemma~\ref{lem:degree-projection} gives the calculation
\begin{equation}\label{eq:Gdegree}
 \deg G_\rho=\sum_i a_i\deg Y_i
 =\sum_i\deg X_i=\delta.
\end{equation}
For the forms defining $\rho$, set
\begin{equation}\label{eq:F}
 F_\rho=G_\rho(\ell_0,\ldots,\ell_{d+1}).
\end{equation}
This is a homogeneous equation of degree $\delta$ for $X$.
It is nonzero because substitution of independent linear forms
induces an injective map of polynomial rings. It vanishes on every
$X_i$, hence on $X$ because $X$ is reduced. Its scalar ambiguity
will not affect either its local order or the linear span of its
derivatives.

For a closed point $x\in\Pj^n_k$, let $\mathfrak n_x$ be the
maximal ideal of $\Ocal_{\Pj^n_k,x}$. If $F$ is homogeneous of
degree $b$, its \emph{order at $x$}, denoted by $\ord_x(F)$, is
the order of the local equation $F/T_i^b$ for any $T_i(x)\ne0$.
Thus $\ord_x(F)\geqslant h$ means
$F/T_i^b\in\mathfrak n_x^h$, and we put $\ord_x(0)=+\infty$.
On changing charts this local equation is multiplied by a unit,
so the definition is independent of the choice. The order is taken
in the ambient local ring; the image of $F_\rho$ in
$\Ocal_{X,x}$ is already zero.

\subsection{Local properties of the projection equations}

To apply the affine Chow-equation statements, let
$\widehat X\subset\A^{n+1}_k$ be the affine cone defined by the
homogeneous ideal of $X$. It is reduced and pure dimensional of
dimension $d+1$. The linear map
\[
 \widehat\rho:\A^{n+1}_k\longrightarrow\A^{d+2}_k,
 \qquad T\longmapsto(\ell_0(T),\ldots,\ell_{d+1}(T))
\]
is allowable for $[\widehat X]$ in the sense of
\cite[Definition~4.1, pp.~324--325]{Kollar1999}. Indeed, the intersection of
the projective closure of $\widehat X$ with the hyperplane at
infinity is $X$, and the center at infinity is precisely
$\Lambda$. The cited criterion identifies allowability with
$\Lambda\cap X=\varnothing$.

The image of the cone over $X_i$ is the cone over $Y_i$.
The degrees of these maps equal $a_i$: choosing a projecting
coordinate nonzero on $X_i$ adjoins the same independent cone
coordinate to $k(X_i)$ and $k(Y_i)$, which does not change their
field-extension degree. Thus the equation of the affine image
cycle is the homogeneous polynomial $G_\rho$, and its pullback is
$F_\rho$. Conversely, every allowable linear projection of the
cone gives an admissible projective projection in this way. In the
hypersurface case the identity gives its defining equation.
Allowing a constant term in the affine projection merely translates
the target and leaves its pulled-back equation unchanged. Consequently
the $F_\rho$ generate the ideal of Chow equations of
$[\widehat X]$ in \cite[Definition~4.1]{Kollar1999}.

We also identify the local multiplicities and orders in this
correspondence. Choose $x\in X$ and a nonzero lift $v\in\widehat X$,
scaled so that $T_i(v)=1$. Writing $t=T_i$ and $u_j=T_j/T_i$ gives
\[
 \widehat X\cap D(T_i)
 \simeq \bigl(X\cap D_+(T_i)\bigr)\times\mathbb G_m.
\]
The point $v$ corresponds to $(x,1)$. Adjoining the local parameter
$t-1$ adjoins one polynomial variable to the associated graded
ring. It follows that
\begin{equation}\label{eq:cone-local-mu}
 \mu_v(\widehat X)=\mu_x(X),
\end{equation}
and $v$ is smooth on $\widehat X$ exactly when $x$ is smooth on
$X$. In the corresponding ambient coordinates,
$F_\rho=t^\delta f_\rho(u)$, where $f_\rho=F_\rho/T_i^\delta$
is the local projective equation. Since $t$ is a unit and
$f_\rho$ is independent of $t$, we also have
\begin{equation}\label{eq:cone-local-order}
 \ord_v(F_\rho)=\ord_x(F_\rho).
\end{equation}
These identifications make no assumption on the characteristic.

\begin{prop}\label{prop:equation-order}
For every closed point $x\in X$ and every admissible projection,
\begin{equation}\label{eq:F-order}
 \ord_x(F_\rho)\geqslant\mu_x(X).
\end{equation}
\end{prop}

\begin{proof}
By \cite[Example~4.6, p.~326]{Kollar1999}, the ideal of Chow equations of
an effective affine cycle is contained in the power of the maximal
ideal specified by its multiplicity at that point. Apply this to
$[\widehat X]$ at $v$ and to its generator $F_\rho$.
Equations~\eqref{eq:cone-local-mu} and
\eqref{eq:cone-local-order} give \eqref{eq:F-order}.
\end{proof}

\begin{prop}\label{prop:smooth-projection}
For every smooth closed point $x\in X$ there is an admissible
projection $\rho$ such that $\ord_x(F_\rho)=1$.
\end{prop}

\begin{proof}
The lift $v$ is smooth on the reduced cone. By
\cite[Example~4.6, p.~326]{Kollar1999}, the ideal of its Chow equations
agrees with its defining ideal near $v$. Since the cone is smooth
of positive codimension there, its defining ideal contains an
element with a nonzero linear term in the ambient local ring.
If every generator $F_\rho$ belonged to the square of the ambient
maximal ideal, their generated ideal would belong to that square
as well. Some $F_\rho$ therefore has order one at $v$, and
\eqref{eq:cone-local-order} gives the assertion at $x$.
\end{proof}

The two propositions supply exactly the local information needed
for the next section. Every projection equation has a prescribed
minimum order along $X$, while at each smooth point at least one
of them has a nonzero first derivative.

\section{A derivative system and the proof of the bound}\label{sec:main-proof}

We now use the projection equations to construct hypersurface
sections of degree $\delta-1$. At each section the degree of the
cycle is multiplied by $\delta-1$, and its local multiplicity at
$x$ is multiplied by at least $\mu_x(X)-1$. We can continue the
construction as long as every component of the current cycle meets
the smooth locus of $X$. The dimension of the singular locus will
therefore determine the permissible number of sections.

\subsection{The derivative system}

We retain the setting of Section~\ref{sec:projection}, with
$1\leqslant d<n$ and $\delta\geqslant2$. The derivatives of a
single projection equation may vanish at a smooth point, because
that particular projection need not preserve smoothness there.
Proposition~\ref{prop:smooth-projection} instead allows the
projection to depend on the point. We consequently take all the
projection equations together and define
\begin{equation}\label{eq:W}
 W=\operatorname{span}_k
 \left\{\frac{\partial F_\rho}{\partial T_j}:
       \rho\text{ admissible},\ 0\leqslant j\leqslant n\right\}.
\end{equation}
Changing the equation of a projected cycle by a scalar multiplies its
partial derivatives by that scalar. Thus $W$ is independent of all
such choices. It is finite dimensional because it belongs to the
degree-$(\delta-1)$ part of $k[T_0,\ldots,T_n]$. In particular, its
common zero set is the zero set of finitely many forms, obtained by
taking a basis of $W$.
Every nonzero member of $W$ defines a hypersurface of degree
$\delta-1$. The following proposition states the local order and
the avoidance property that will be used to choose these
hypersurfaces successively.

\begin{prop}\label{prop:W}
The system $W$ has the following properties.
\begin{enumerate}[(i)]
\item At every closed point $x\in X$, every $g\in W$ has local order at least
$\mu_x(X)-1$.
\item At every smooth closed point of $X$, some $g\in W$ is nonzero.
In particular, the base locus of $W$ on $X$, defined by the common
zero set of its members, is contained in $\Sing(X)$.
\end{enumerate}
\end{prop}

\begin{proof}
For every admissible projection, Proposition~\ref{prop:equation-order}
and Lemma~\ref{lem:order} give
$\mu_x(V(F_\rho))=\ord_x(F_\rho)\geqslant\mu_x(X)$.
If $\mu_x(X)>1$, the lower-bound assertion of
\cite[Corollary~5.5, p.~348]{Liu2022}, applied with derivative order
one to the hypersurface $V(F_\rho)$, shows that every nonzero
first derivative has local order at least $\mu_x(X)-1$.
For $\mu_x(X)=1$ the required lower bound is zero.
The corresponding power of the ambient maximal ideal is closed
under linear combinations, so the bound holds for every $g\in W$.
This proves (i).

At a smooth point $x$, Proposition~\ref{prop:smooth-projection}
gives an $F_\rho$ of local order one.
By \cite[Corollary~5.4, p.~347]{Liu2022}, some first derivative
of $F_\rho$ is nonzero at $x$. It belongs to $W$, proving (ii).
The assertion about the base locus follows because every nonempty
locally closed subset of a scheme of finite type over the
algebraically closed field $k$ contains a closed $k$-point.
\end{proof}

These first-derivative properties hold in every characteristic.
In the chart $T_0\ne0$, with
$f_\rho(t)=F_\rho(1,t_1,\ldots,t_n)$, the derivatives with index
$j\geqslant1$ dehomogenize to $\partial f_\rho/\partial t_j$.
The remaining derivative is related to them by
\begin{equation}\label{eq:Euler}
 \left(\frac{\partial F_\rho}{\partial T_0}\right)(1,t)
 =\delta f_\rho(t)-\sum_{j=1}^n
       t_j\frac{\partial f_\rho}{\partial t_j}(t).
\end{equation}
This identity involves no division by $\delta$, so no assumption
on its residue in $k$ is needed.

\subsection{Successive hypersurface sections}\label{subsec:proper-sections}

We work over an algebraically closed field $k$, in the setting of
Proposition~\ref{prop:W}. Thus $X$ is reduced of pure dimension
$1\leqslant d<n$, its degree is $\delta\geqslant2$, and $W$ is the
space of first derivatives in \eqref{eq:W}. We shall use the two
properties of $W$ at different stages. Its base locus determines how
many proper sections can be taken. The order of its members determines
the multiplicities of the resulting cycle.

\begin{prop}\label{prop:aux-cycle}
Let $1\leqslant t\leqslant d$ satisfy $\dim\Sing(X)\leqslant d-t$.
There are $g_1,\ldots,g_t\in W$ such that, putting $Z_0=[X]$ and
$Z_j=Z_{j-1}\cdot V(g_j)$ for $1\leqslant j\leqslant t$, each
intersection is proper on the components of $Z_{j-1}$. For every
$j=1,\ldots,t$, $Z_j$ is effective of pure dimension $d-j$ and
\begin{equation}\label{eq:Zdegree}
 \deg Z_j=\delta(\delta-1)^j.
\end{equation}
At every closed point $x$ of $X$ one has
\begin{equation}\label{eq:Zmu}
 \mu_x(Z_j)\geqslant\mu_x(X)\bigl(\mu_x(X)-1\bigr)^j .
\end{equation}
\end{prop}

\begin{proof}
Suppose that $Z_j$ has been constructed, with $0\leqslant j<t$.
An irreducible component $V$ of $Z_j$ has dimension
$d-j\geqslant d-t+1$. Thus $V$ meets the smooth locus of $X$.
The intersection is open in $V$ and contains a closed point $z$,
since $k$ is algebraically closed. By Proposition~\ref{prop:W} some
member of $W$ does not vanish at $z$. Hence the subspace
\[
 W_V=\{g\in W\mid g|_V=0\}
\]
is proper. The point $z$ and the member exhibiting properness may
depend on $V$.

There are only finitely many components of $Z_j$. For each, choose
a nonzero linear functional on $W$ that vanishes on $W_V$. Their
product is a nonzero polynomial on $W$. Since $k$ is infinite, we
can choose $g_{j+1}\in W$ on which this product does not vanish.
It cuts all the components properly in one hypersurface. Its degree
is $\delta-1$, so Lemma~\ref{lem:bezout} gives an effective pure
cycle of dimension $d-j-1$ and degree $(\delta-1)\deg Z_j$.
As $\deg Z_0=\delta$, this proves \eqref{eq:Zdegree}.

Fix a closed point $x\in X$, and put $u=\mu_x(X)$. If $u=1$, the
right hand side of \eqref{eq:Zmu} is zero for each $j\geqslant1$.
If $u\geqslant2$, every $g_{j+1}$ has ambient order at least $u-1$
at $x$, by Proposition~\ref{prop:W}. Write
$Z_j=\sum_i b_i[V_i]$. For each $V_i$ through $x$, the quotient map
from the ambient local ring to $\Ocal_{V_i,x}$ sends the ambient
maximal ideal to the maximal ideal of $\Ocal_{V_i,x}$. Hence the
local equation of $V(g_{j+1})$ belongs to the $(u-1)$st power of
this maximal ideal. It is nonzero on $V_i$ because the section was
chosen properly. Corollary~\ref{cor:cycle-cut} therefore gives
\[
 \begin{split}
 \mu_x(Z_{j+1})
 &=\sum_{i:x\in V_i} b_i\,
      \mu_x\bigl([V_i]\cdot V(g_{j+1})\bigr)\\
 &\geqslant\sum_{i:x\in V_i}b_i(u-1)\mu_x(V_i)
  =(u-1)\mu_x(Z_j).
 \end{split}
\]
Starting with $\mu_x(Z_0)=u$, the first section gives $u(u-1)$,
and each further section contributes another factor $u-1$.
This proves \eqref{eq:Zmu}. Since $u\geqslant2$, this lower bound
is positive, so $x$ remains in the support at each step. Forming
the intersections component by component accounts
for the generic lengths in $Z_j$, even if an intermediate
scheme-theoretic intersection has embedded components.
\end{proof}

For $t=d-s$, where $s=\dim\Sing(X)\geqslant0$, the last cycle has
dimension $s$. At this stage an entire component can lie in the
singular locus, and the preceding argument no longer supplies a
proper section on that component. This is where the codimension of
the singular locus enters the bound.

\subsection{Extension of the ground field}\label{subsec:basechange}

Return to a scheme $X$ over $\F_q$ and let $k=\kbar$. For a rational
point $x\colon\Spec\F_q\longrightarrow X$, its extension of scalars
is a single $k$-point $\bar x$ of $\overline X=X\times_{\F_q}k$.
We shall take the sections over $k$ and sum their multiplicities
over just these points. There are two issues to check. The
geometric construction must preserve the degree, dimension, and
local multiplicities occurring in the original problem. Moreover,
the point-counting estimate must still involve $q$, even though
the auxiliary cycle is defined over $k$. The last assertion of the
following proposition makes this distinction explicit.

\begin{prop}\label{prop:finite-field-transfer}
Let $X\subset\Pj^n_{\F_q}$ be a nonempty reduced closed subscheme of
pure dimension $d$ and degree $\delta$. Then $\overline X$ is reduced
of pure dimension $d$ and degree $\delta$, and
$\dim\Sing(\overline X)=\dim\Sing(X)$. Moreover
\begin{equation}\label{eq:basechange-mu}
 \mu_{\bar x}(\overline X)=\mu_x(X)\qquad(x\in X(\F_q)).
\end{equation}
Suppose that $t\geqslant1$ and that $C$ is an effective pure $u$-cycle
in $\Pj^n_k$ such that, for every $x\in X(\F_q)$,
\[
 \mu_{\bar x}(C)\geqslant
 \mu_{\bar x}(\overline X)\bigl(\mu_{\bar x}(\overline X)-1\bigr)^t.
\]
Then
\begin{equation}\label{eq:transfer}
 S_t(X;q)\leqslant(\deg C)\pi_u(q).
\end{equation}
\end{prop}

\begin{proof}
Since $\F_q$ is perfect, a reduced scheme of finite type over
$\F_q$ is geometrically reduced \cite[Corollaire~4.6.11]{EGAIV2}.
The geometric components lying over a given component of $X$ all have
its dimension \cite[Proposition~4.2.7(iv)]{EGAIV2}. In particular every
component of $\overline X$ has dimension $d$, even if an irreducible
component of $X$ splits. If $R$ is the homogeneous coordinate ring,
the dimension of $(R\otimes_{\F_q}k)_N$ over $k$ equals that of $R_N$
over $\F_q$ in every degree $N$. The Hilbert polynomial is unchanged
and therefore so is its leading coefficient $\delta/d!$.

The smooth locus pulls back to the smooth locus
\cite[Corollaire~6.7.8(iv) and Th\'eor\`eme~6.8.6]{EGAIV2}.
Consequently the inverse image of the singular locus has underlying
closed set $\Sing(\overline X)$. Applying invariance of dimension
under field extension to the reduced closed subscheme $\Sing(X)_{\rm red}$
gives the asserted dimension equality \cite[Corollaire~4.1.4]{EGAIV2}.
It holds for the empty locus by our convention on its dimension.

The invariance of multiplicity under a finite extension of
$\F_q$ is \cite[Proposition~2.6, pp.~325--326]{Liu2022}, since every such
extension is Galois. The passage to the algebraic closure is
compatible with the same Hilbert--Samuel functions. Indeed, for
an affine neighborhood $\Spec A$ of a rational point $x$, with
maximal ideal $\m$, and any extension $K/\F_q$, put
$A_K=A\otimes_{\F_q}K$ and $\m_K=\m A_K$. Since $A/\m=\F_q$,
we have $A_K/\m_K=K$ and
\begin{equation}\label{eq:basechange-length}
 A_K/\m_K^{N+1}\simeq(A/\m^{N+1})\otimes_{\F_q}K.
\end{equation}
The dimensions of these quotients over their residue fields,
and hence their local lengths, agree. Taking $K=k$ gives
\eqref{eq:basechange-mu}.

The assignment $x\mapsto\bar x$ identifies $X(\F_q)$ with a subset
of the ambient geometric points of $\overline X$. By the hypothesis
on $C$, \eqref{eq:basechange-mu}, and nonnegativity of the other
multiplicities,
\[
 S_t(X;q)\leqslant\sum_{x\in X(\F_q)}\mu_{\bar x}(C)
 \leqslant\sum_{x\in\Pj^n(\F_q)}\mu_{\bar x}(C).
\]
Proposition~\ref{prop:projective-count} applies to this last sum
over ambient rational points even if a component of $C$ is not
defined over $\F_q$. It gives \eqref{eq:transfer}.
\end{proof}

\subsection{Proof of the main theorem and consequences}

The auxiliary cycle now supplies the connection between the local
weights and the global point count. Its dimension determines the
power of $q$ in the estimate, while its degree determines the factor
involving $\delta$.

\begin{proof}[Proof of Theorem~\ref{thm:main}]
If $X$ is smooth, every rational point has multiplicity one and
$S_t(X;q)=0$. Otherwise Proposition~\ref{prop:finite-field-transfer}
gives a reduced scheme $\overline X$ of pure dimension $d$ and degree
$\delta$ with a nonempty singular locus of dimension at most $d-t$.
Here $d<n$, since a reduced closed subscheme of dimension $n$ in
$\Pj^n_k$ is $\Pj^n_k$ itself. Also $\delta\geqslant2$ by
Lemma~\ref{lem:degree-one}. The hypotheses for the construction of
$W$ are therefore satisfied.

Apply Proposition~\ref{prop:aux-cycle} to $\overline X$. At every
$\bar x$ above a rational point, \eqref{eq:Zmu} bounds the
summand of $S_t(X;q)$ by the multiplicity of $Z_t$. Since this cycle
has pure dimension $d-t$, Proposition~\ref{prop:finite-field-transfer}
and \eqref{eq:Zdegree} give
\[
 S_t(X;q)\leqslant(\deg Z_t)\pi_{d-t}(q)
       =\delta(\delta-1)^t\pi_{d-t}(q).
\]
\end{proof}

\subsubsection{}\label{subsec:single-section}
For the quadratic sum, the construction has a direct interpretation
without singular-locus data. Every component of the reduced
$\overline X$ meets its smooth locus. Hence one member $g$ of $W$
cuts all components of $[\overline X]$ properly. If a geometric closed point $x$ has
multiplicity $u\geqslant2$, the local equation of this single
section has order at least $u-1$. Corollary~\ref{cor:cycle-cut}
gives
\[
 \mu_x([\overline X]\cdot V(g))\geqslant u(u-1),
 \qquad
 \deg([\overline X]\cdot V(g))=\delta(\delta-1).
\]
Applying the projective count in dimension $d-1$ gives
Corollary~\ref{cor:reduced}. This calculation uses density of the
smooth locus only on the original components. For subsequent
sections the corresponding density on the current components follows
from their dimension, as in Proposition~\ref{prop:aux-cycle}.

\subsubsection{}\label{subsec:field-size}
A parameter space with a nonempty open subset over $\F_q$ can have
no rational point in that open subset. For example,
$\A^1_{\F_q}\setminus V(T^q-T)$ has this property. The choice of
projection and the simultaneous avoidance of the subspaces $W_V$
were both made over $k$. Although the finitely many coefficients of
the chosen sections belong to some common finite extension of
$\F_q$, its cardinality never occurs in the count. In
Proposition~\ref{prop:affine-count} the coordinate of an ambient
$\F_q$-point still takes only $q$ possible values. In particular the
argument also applies for $\delta>q$. When the characteristic divides
$\delta$, the same proof works because the field degrees in the image
cycle include inseparable degrees and the homogeneous derivative
identity \eqref{eq:Euler} requires no division by $\delta$.

\begin{rema}\label{rem:hypersurface-proof}
For a reduced hypersurface with equation $F$, the identity projection
gives $F$ itself. On each affine chart the order of the
dehomogenized equation equals the local multiplicity by
Lemma~\ref{lem:order}. A first derivative consequently has order at
least $\mu_x(X)-1$ and at a smooth point some derivative has order
zero. Thus one can take $W$ to be the span of the partial derivatives
of this single $F$ throughout the proof. It is the weighted count
of the proper intersection cycle that gives \eqref{eq:hypersurface}
in place of the sum of terms in \eqref{eq:old}.
\end{rema}

\begin{coro}\label{cor:many}
Under the hypotheses of Theorem~\ref{thm:main}, for every integer
$M\geqslant2$ one has
\begin{equation}\label{eq:large-mult}
 \#\{x\in X(\F_q)\mid\mu_x(X)\geqslant M\}
 \leqslant \frac{\delta(\delta-1)^t\pi_{d-t}(q)}{M(M-1)^t}.
\end{equation}
Without a hypothesis on the dimension of the singular locus, the
same count is at most $\delta(\delta-1)\pi_{d-1}(q)/(M(M-1))$.
If the singular locus is nonempty of dimension $s$, put $r=d-s$.
For every $1\leqslant u\leqslant r$ one has
\begin{equation}\label{eq:lower-moments}
 S_u(X;q)\leqslant
 \min_{u\leqslant t\leqslant r}
 \bigl\{\delta(\delta-1)^t\pi_{d-t}(q)\bigr\}.
\end{equation}
\end{coro}

\begin{proof}
A point counted on the left of \eqref{eq:large-mult} contributes at
least $M(M-1)^t$ to $S_t(X;q)$. Dividing the bound in
Theorem~\ref{thm:main} by that number proves the assertion. The
second estimate follows by taking $t=1$ in Corollary~\ref{cor:reduced}.

For the last assertion fix $t$ between $u$ and $r$. At a point
of multiplicity one both summands vanish. At a point of multiplicity
$v\geqslant2$ we have $v-1\geqslant1$, so
$v(v-1)^u\leqslant v(v-1)^t$. Hence $S_u(X;q)\leqslant S_t(X;q)$.
Since $t\leqslant r$, we have $s\leqslant d-t$, and
Theorem~\ref{thm:main} bounds $S_t(X;q)$. This holds for each
integer in the displayed minimum.
\end{proof}

\section{Equality cases and examples beyond hypersurfaces}\label{sec:sharpness}

In \cite[Example~5.11]{Liu2022}, the singular locus of a
hypersurface is a linear space of codimension one in the
hypersurface, and every rational point of that space has
multiplicity equal to the degree. To obtain equality in the present
bound, we use cones whose vertex is a projective space of prescribed
dimension. Their smooth points contribute zero to the sum, while
every point of the vertex has the same multiplicity. The degree and
the number of rational points of the vertex then account separately
for the two factors in the bound. We give the local calculation of
this multiplicity, followed by explicit bases in every characteristic
and an example in codimension two. Before these constructions, we
explain the role of pure dimensionality.

\begin{rema}[Why pure dimensionality is retained]\label{rem:pure-needed}
The pure-dimensional hypothesis in Corollary~\ref{cor:reduced}
is necessary when the degree is defined by the leading coefficient
of the Hilbert polynomial. Fix $d\geqslant2$ and use homogeneous
coordinates
$A_0,\ldots,A_d,B_0,B_1,B_2$ on $\Pj^{d+3}_{\F_q}$.
The subspaces
\[
 L=V(B_0,B_1,B_2)\simeq\Pj^d,
 \qquad \Pi=V(A_0,\ldots,A_d)\simeq\Pj^2
\]
are disjoint. In $\Pi$, let $C=V(B_0B_1)$ be the union of its two
distinct coordinate lines, and set $X=L\cup C$ with the reduced
scheme structure. It is reduced of dimension $d$ but is not pure.
We compute its degree. If $I_L,I_C$
are the two homogeneous ideals in the ambient polynomial ring $S$,
the exact graded sequence
\[
 0\longrightarrow S/(I_L\cap I_C)
 \longrightarrow S/I_L\oplus S/I_C
 \longrightarrow S/(I_L+I_C)\longrightarrow0
\]
has last term $\F_q$, because $I_L+I_C$ contains all the coordinates.
The last term has zero homogeneous components in positive degrees.
It follows that the Hilbert polynomial of $X$ is the sum of those
of $L$ and $C$. The latter has degree one, less than $d$.
The coefficient of the degree-$d$ term is therefore $1/d!$, the
leading coefficient for $L$. In particular, $\deg X=1$.

The point $P=[0:\cdots:0:0:0:1]$ is the intersection of the two
lines. On the chart $B_2\ne0$ of $\Pi$, its local equation is
$uv=0$ in a smooth plane. Lemma~\ref{lem:order} gives
$\mu_P(C)=2$. Since $L$ is disjoint from $C$, the local rings of
$X$ and $C$ at $P$ agree. Every other point is smooth on its
component, so
\[
 \sum_{x\in X(\F_q)}\mu_x(X)(\mu_x(X)-1)=2,
 \qquad \delta(\delta-1)\pi_{d-1}(q)=0.
\]
The two lines through $P$ affect the multiplicity sum, although
neither contributes to the leading coefficient defining $\deg X$.
\end{rema}

\subsection{Cones over a smooth base}

Let the vertex have dimension $s$, and let its codimension in the
cone be $r$. The base then has dimension $r-1$. This notation makes
$r$ the exponent of $\mu_x(X)-1$ in the sharp bound. The same
construction applies to a smooth base of any codimension, and the
base need not be geometrically connected.

\begin{prop}\label{prop:cones}
Let $Y\subset\Pj^M_{\F_q}$ be a nonempty smooth projective scheme
of pure dimension $r-1$, where $r\geqslant1$, and degree $\delta>1$.
Let $R$ be its standard graded homogeneous coordinate ring. For
$s\geqslant0$, define
\[
 X=\Proj R[U_0,\ldots,U_s]\subset\Pj^{M+s+1}_{\F_q},
\]
where each $U_j$ has degree one. This is the projective cone over
$Y$ with vertex
$L=(T_0=\cdots=T_M=0)\simeq\Pj^s_{\F_q}$, where $T_i$ are the
coordinates on $\Pj^M$. The cone is geometrically reduced and pure,
and satisfies
\[
 \dim X=r+s,\qquad \deg X=\delta,\qquad
 \Sing(X)_{\mathrm{red}}=L,
 \qquad\mu_x(X)=\delta\quad(x\in L(\F_q)).
\]
Consequently equality holds in \eqref{eq:main}, and in
\eqref{eq:family} with $t=r$. For $r=1$, equality also holds in
\eqref{eq:quadratic}.
\end{prop}

\begin{proof}
Write $R=\F_q[T_0,\ldots,T_M]/I$, where $I$ is the saturated
homogeneous ideal of $Y$. We first verify that the homogeneous
coordinate ring of the cone remains reduced over $k=\kbar$. The
map
\[
 R\longrightarrow\prod_{i=0}^M R[T_i^{-1}]
\]
is injective. In fact, an element in its kernel is killed by a power
of each $T_i$, and hence by a power of the ideal
$(T_0,\ldots,T_M)$. The saturation of $I$ gives the assertion.
Tensor with $k$. Since $k/\F_q$ is flat and the product is finite,
the resulting map is still injective. Each nonzero factor is a Laurent
polynomial ring in $T_i$ over the coordinate ring of the chart
$Y_k\cap D_+(T_i)$. Indeed, a homogeneous fraction of degree
$m$ in the localization is the product of a fraction of degree zero
and $T_i^m$. The coordinate ring of this affine chart is reduced since
$Y_k$ is smooth. The product on the right is therefore reduced, so
$R_k=R\otimes_{\F_q}k$ and its polynomial ring
$R_k[U_0,\ldots,U_s]$ are reduced.

Every minimal prime of $R_k[U_0,\ldots,U_s]$ is extended from a
minimal prime of $R_k$. To see this, write the zero ideal of $R_k$
as the finite intersection of its minimal primes. Extension to a
polynomial ring preserves this intersection coefficient by coefficient,
and each extended prime is again prime. Thus the geometric components
of $X$ are the cones over the components of $Y_k$. Smoothness and
purity give dimension $r-1$ for each component of $Y_k$. To compute
the dimension and degree of the cone, we use its Hilbert function.

Adding one variable of degree one replaces a Hilbert function $h(N)$
by $\sum_{j=0}^N h(j)$. If $h(N)$ has leading term
$\delta N^{r-1}/(r-1)!$, its cumulative sum has leading term
$\delta N^r/r!$. This follows by summing powers, or by taking the
finite difference of the claimed leading polynomial. Repeating
$s+1$ times shows that $R[U_0,\ldots,U_s]$ has Hilbert polynomial
of degree $r+s$ and leading coefficient $\delta/(r+s)!$.
Thus $\dim X=r+s$ and $\deg X=\delta$. This also covers $r=1$,
when the Hilbert polynomial of $R$ is the constant $\delta$.

Applying the same calculation to each minimal-prime quotient gives
dimension $r+s$ for every geometric component of $X$. To locate its
singular points, consider first the chart $T_i\ne0$ on $X_k$. Its
coordinate ring is
\[
 (R_k[T_i^{-1}])_0
 [U_0/T_i,\ldots,U_s/T_i].
\]
This chart is a product of a smooth affine chart of $Y_k$ and
$\A^{s+1}_k$. In particular, the additional coordinates of the cone
introduce no singularities on this chart. Polynomial algebras are
smooth, and smoothness is stable under base change and composition
\cite[Chapter III, Proposition 10.1(b)--(c)]{Hartshorne}.
These charts cover $X_k\setminus L_k$.

Now choose any closed point $x\in L_k$ and an index $j$ with
$U_j(x)\ne0$. On the chart $U_j\ne0$, dehomogenization identifies
the coordinate ring with
\[
 R_k[v_1,\ldots,v_s].
\]
Indeed, the homogeneous relations among the original coordinates
$T_i$ become the same relations among $T_i/U_j$. The remaining
ratios $U_h/U_j$ are independent polynomial variables, denoted here
by $v_1,\ldots,v_s$. All the $T_i/U_j$ vanish at $x$. Since $k$ is
algebraically closed, the values $v_h(x)$ belong to $k$. Subtracting
these values, write $w_h=v_h-v_h(x)$. The local ring is the
localization of the standard graded ring
\[
 S=R_k[w_1,\ldots,w_s]
\]
at its homogeneous maximal ideal $S_+$. As $S$ is generated in degree
one, $S_+^a=\bigoplus_{b\geqslant a}S_b$. The quotient $S/S_+^{N+1}$
has residue field $k$ and nilpotent positive-degree ideal. It is
therefore already local Artinian. Localization at $S_+$ has no
effect, and its length is its $k$-dimension. Thus
\begin{equation}\label{eq:cone-HS}
 \length_{S_{S_+}}(S_{S_+}/S_+^{N+1}S_{S_+})
 =\sum_{a=0}^N\dim_k S_a.
\end{equation}
The sum on the right is the dimension of the degree-$N$ part of
$S[z]$ with $\deg z=1$. As graded rings, $S[z]$ and
$R_k[U_0,\ldots,U_s]$ are isomorphic. The computation of its
Hilbert polynomial above therefore gives leading term
$\delta N^{r+s}/(r+s)!$ in \eqref{eq:cone-HS}. By the definition
of Hilbert--Samuel multiplicity, $\mu_x(X_k)=\delta$.
Since $\delta>1$ and regular local rings have multiplicity one,
every such $x$ is singular. Together with smoothness off $L_k$, this
proves equality of the underlying singular locus with $L_k$; it
descends to the stated equality of reduced closed loci over $\F_q$.

For rational points on $L$, multiplicity is $\delta$ by
\eqref{eq:basechange-mu}. Off $L$ it is one. Therefore, for every
integer $t\geqslant1$,
\[
 S_t(X;q)=\sum_{x\in L(\F_q)}\delta(\delta-1)^t
          =\delta(\delta-1)^t\pi_s(q).
\]
Taking $t=r$ gives $d-t=s$, since $d=\dim X=r+s$. The last
expression is then exactly the right-hand side of
\eqref{eq:family}, as required.
\end{proof}

The calculation does not require the base to have any rational
points. Its degree enters through the Hilbert--Samuel multiplicity
at the vertex, whereas all rational points contributing to the sum
lie on the vertex itself. Thus a smooth base of degree $\delta$
suffices even when $\delta$ is large compared with $q$.

\subsection{Hypersurfaces of prescribed degree}

\begin{exem}[Sharpness within hypersurfaces]\label{ex:fermat}
Fix a finite field $\F_q$ of characteristic $p$ and integers
$r\geqslant1$, $s\geqslant0$ and $\delta\geqslant2$.
We construct a smooth hypersurface $Y\subset\Pj^r_{\F_q}$ of degree
$\delta$. Proposition~\ref{prop:cones} then gives a hypersurface of
dimension $r+s$ attaining \eqref{eq:family} with $t=r$, for these
prescribed values of $q,r,s$ and $\delta$. The usual Fermat
equation treats degrees prime to the characteristic. We give separate
equations for the remaining degrees so that the equality examples
impose no condition on the characteristic or the size of the field.

For the smoothness computations, let $F$ be a homogeneous equation
and consider a geometric zero of $F$ in the chart $T_j\ne0$. After
normalizing $T_j$ to one, the affine partial derivatives are the
homogeneous partials with indices $i\ne j$ evaluated at this point.
Euler's identity $\sum_i T_i\partial F/\partial T_i=\delta F$ shows
that, if all these partials vanish, the remaining partial also
vanishes. This argument holds even when $p\mid\delta$. Thus it
suffices to show that no projective zero of $F$ annihilates all the
homogeneous partial derivatives. Indeed, a nonzero affine partial
means that the local equation has nonzero class in the cotangent
space of the ambient regular local ring. It extends to a regular
system of parameters, so the hypersurface is regular at that point
\cite[Theorem~14.2]{Matsumura}. This verifies smoothness at all
geometric closed points, which suffices since the nonsmooth locus
is closed. For each equation below, multiplication by $F$ on the
polynomial ring in $r+1$ variables is injective. Taking the
difference of its Hilbert polynomials in degrees $m$ and
$m-\delta$ gives dimension $r-1$ and degree $\delta$ for $Y$.

If $p\nmid\delta$, take
\[
 F=T_0^\delta+\cdots+T_r^\delta.
\]
Its partial derivatives are $\delta T_i^{\delta-1}$; their
simultaneous vanishing would force every coordinate to vanish.
This gives the desired smooth base in degrees prime to $p$.

Suppose that $p\mid\delta$ and $\delta>2$. Take
\[
 F=T_0^\delta+\sum_{i=1}^r T_{i-1}T_i^{\delta-1}.
\]
The first partial derivative is $T_1^{\delta-1}$. At any common
zero of the partial derivatives we therefore have $T_1=0$. For
$1\leqslant i<r$ the derivative with respect to $T_i$ is
\[
 (\delta-1)T_{i-1}T_i^{\delta-2}+T_{i+1}^{\delta-1}.
\]
As $\delta-2>0$, if $T_i=0$ then vanishing of this derivative
gives $T_{i+1}=0$. Induction gives $T_1=\cdots=T_r=0$. The
remaining projective point is $[1:0:\cdots:0]$, at which $F=1$,
so it does not lie on $Y$.

It remains to treat $p=\delta=2$. Pair consecutive coordinates
and put
\[
 F=\begin{cases}
 \displaystyle\sum_{j=0}^{(r-1)/2}T_{2j}T_{2j+1},&r\text{ odd},\\[6pt]
 \displaystyle\sum_{j=0}^{r/2-1}T_{2j}T_{2j+1}+T_r^2,&r\text{ even}.
 \end{cases}
\]
The derivative with respect to a paired coordinate is its partner. Vanishing of
all the partials thus forces all the paired coordinates to vanish.
If $r$ is odd, there is then no projective point left. If $r$
is even, the additional equation $F=0$ gives $T_r^2=0$ and again
there is no projective point. This completes the construction in
every characteristic.
\end{exem}

\begin{exem}[The quadratic bound]\label{ex:quadratic-sharp}
For the case $t=1$, one can also see the behavior when $p\mid\delta$
on a binary form. A smooth zero-dimensional subscheme of degree
$\delta$ in $\Pj^1$ is a suitable base. When $p\nmid\delta$ take
$T_0^\delta+T_1^\delta$. When $p\mid\delta$ take
\[
 T_0^\delta-T_0T_1^{\delta-1}.
\]
Its zero scheme is contained in the affine chart $T_1\ne0$;
in that chart the derivative of $z^\delta-z$ is $-1$. Thus all
its $\delta$ geometric zeros are simple. For $\delta=p$ the
factorization
\[
 T_0^p-T_0T_1^{p-1}=\prod_{a\in\F_p}(T_0-aT_1)
\]
shows that these $p$ points are already rational.

Use the same binary equation in $\Pj^{d+1}_{\F_q}$, with
$d\geqslant1$. Its zero scheme is a cone with vertex
$L=V(T_0,T_1)\simeq\Pj^{d-1}_{\F_q}$. Proposition~\ref{prop:cones}
gives multiplicity $\delta$ at every rational point of $L$ and
multiplicity one off $L$. In particular,
\[
 S_1(X;q)=\delta(\delta-1)\#L(\F_q)
          =\delta(\delta-1)\pi_{d-1}(q).
\]
For each fixed $q$ and $d$, equality in \eqref{eq:quadratic}
therefore holds in every degree $\delta\geqslant2$. The degrees
$\delta=mp+1$, for example, can be arbitrarily large compared to
$q$.
\end{exem}

\subsection{A codimension-two example which is not a complete intersection}

The preceding equality examples are hypersurfaces. We now use a
smooth curve of codimension two as base to show that equality also
occurs for schemes which are not local complete intersections.
Take the twisted cubic. We describe its coordinate ring explicitly,
both to compute its degree and to determine the number of equations
needed at a vertex point of its cone. Consider
\[
 \nu_3:\Pj^1\longrightarrow\Pj^3,
 \qquad [u:v]\longmapsto[u^3:u^2v:uv^2:v^3],
\]
and write the four target coordinates as $a,b,c,d$. Its image $Y$ has
homogeneous ideal
\begin{equation}\label{eq:twisted-ideal}
 I=(b^2-ac,\ bc-ad,\ c^2-bd).
\end{equation}
We verify this over an arbitrary field $k$, including characteristics
two and three. Give the variables weights $(0,1,1,0)$. In the three
displayed relations the monomials $b^2$, $bc$, and $c^2$ can
respectively be replaced by $ac$, $ad$, and $bd$;
each replacement strictly lowers the weight and preserves total
degree. Repeated replacements therefore terminate. Every class of
degree $m$ is spanned by monomials of the following forms:
\[
 a^i d^j\ (i+j=m),\qquad
 a^i b d^j\ (i+j=m-1),\qquad
 a^i c d^j\ (i+j=m-1).
\]
There are $(m+1)+m+m=3m+1$ of them for $m\geqslant1$, and one for
$m=0$. Under substitution by the four cubic monomials above, their
$v$-exponents are respectively $3j$, $3j+1$, and $3j+2$.
For the first family $0\leqslant j\leqslant m$, whereas for the
other two $0\leqslant j\leqslant m-1$. The exponents are therefore
all distinct and together range from $0$ to $3m$. Their
images are thus a basis of the binary forms of degree $3m$.
It follows degree by degree that the substitution map induces an
isomorphism
\begin{equation}\label{eq:veronese}
 k[a,b,c,d]/I\simeq\bigoplus_{m\geqslant0}k[u,v]_{3m},
\end{equation}
where the right side is graded by $m$. In particular $I$ is prime,
is the full homogeneous ideal of the image, and the Hilbert polynomial
is $3m+1$. Thus $\dim Y=1$ and $\deg Y=3$.

The charts $a\ne0$ and $d\ne0$ cover $Y$. For if $a=d=0$, the
relations give $b^2=c^2=0$, and hence $b=c=0$ at every geometric
point, which is impossible in projective space. On $a\ne0$, putting
$t=b/a$ gives $c/a=t^2$ and $d/a=t^3$, so this chart has ring $k[t]$.
On $d\ne0$, set $w=c/d$. The equations give $b/d=w^2$ and
$a/d=w^3$, hence this chart has coordinate ring $k[w]$. On the
overlap, $tw=bc/(ad)=1$ by the middle equation in
\eqref{eq:twisted-ideal}. Thus
$Y\simeq\Pj^1$, in particular it is smooth over every field.

We now take $k=\F_q$, adjoin $s+1$ variables
$U_0,\ldots,U_s$, and form the projective cone
$X\subset\Pj^{s+4}_{\F_q}$ defined by the same three quadrics.
Proposition~\ref{prop:cones} gives
\begin{equation}\label{eq:twisted-count}
 \begin{gathered}
 \dim X=s+2,\qquad\deg X=3,\qquad\Sing(X)_{\mathrm{red}}=\Pj^s,\\
 S_2(X;q)=12\pi_s(q),\qquad S_1(X;q)=6\pi_s(q).
 \end{gathered}
\end{equation}
The ambient space has dimension $s+4$, so $X$ has codimension two.
We check that it is not a local complete intersection at a vertex
point. Use the chart $U_0\ne0$ and normalize $U_0$ to one. We keep
the letters $a,b,c,d,U_1,\ldots,U_s$ for the resulting affine
coordinates, and take the vertex point at which all of them vanish.
The ambient regular local ring has these $s+4$ coordinates as a
regular system of parameters. Denote its maximal ideal by
$\mathfrak n$. The defining ideal $J$ is
generated by the three quadrics in \eqref{eq:twisted-ideal} and is
contained in $\mathfrak n^2$. Their classes in $J/\mathfrak nJ$ are
linearly independent. In fact, $\mathfrak nJ\subseteq\mathfrak n^3$,
whereas each of the monomials $b^2$, $bc$, $c^2$ occurs in just one
of the three degree-two initial forms. A linear combination of the
quadrics with coefficients in $k$ can therefore belong to
$\mathfrak nJ$ only when all three coefficients vanish. The classes
also span $J/\mathfrak nJ$, since the quadrics generate $J$.
Consequently
\[
 \dim_k J/\mathfrak nJ=3.
\]
Nakayama's lemma shows that $J$ requires three generators. Its
height is two, by the dimension $s+2$ computed above. A local
complete intersection of codimension two in this regular ambient
space would have a defining regular sequence of length two.
This would give two generators of $J$, contrary to the computation
of $J/\mathfrak nJ$. Taking $s=0$ gives a surface in $\Pj^4$
with a single singular point of multiplicity three. At that point
the weight in $S_2$ is $3(3-1)^2=12$, which accounts for the entire
sum.

\subsection{The best constant when the ambient dimension is fixed}

The normalization in \cite[Conjecture~5.13]{Liu2022} divides the
sum by $\delta^{d-s+1}q^s$. For fixed ambient dimension, the cone
calculation determines the least constant in this normalization
as well.

\begin{coro}\label{cor:optimal-n}
For $n\geqslant2$, the smallest real number $C_n$ such that
\begin{equation}\label{eq:Cn}
 \sum_{x\in X(\F_q)}\mu_x(X)(\mu_x(X)-1)^{d-s}
 \leqslant C_n\delta^{d-s+1}q^s
\end{equation}
for every prime power $q$ and every reduced pure-dimensional
$X\subset\Pj^n_{\F_q}$ of dimension $d$, degree $\delta$ and
nonempty singular locus of dimension $s$ is
\[
 C_n=2-2^{-(n-2)}.
\]
\end{coro}

\begin{proof}
A reduced closed subscheme of dimension $n$ in $\Pj^n$
is $\Pj^n$ itself. As $X$ has singular points, it follows that
$d\leqslant n-1$. The smooth locus contains the generic points of
$X$, whence $s\leqslant d-1\leqslant n-2$. Corollary~\ref{cor:conjecture} therefore gives
\[
 \frac{\sum_x\mu_x(X)(\mu_x(X)-1)^{d-s}}
      {\delta^{d-s+1}q^s}
 \leqslant 2-2^{-s}\leqslant2-2^{-(n-2)}.
\]
For the reverse inequality, work over $\F_2$, put $s=n-2$ and
$r=1$, and let $\delta$ run through odd integers greater than one.
The polynomial $T_0^\delta+T_1^\delta$ defines a smooth
zero-dimensional scheme of degree $\delta$ in $\Pj^1$, since its
two derivatives have no common projective zero. Its cone with vertex
$\Pj^{n-2}$ is a reduced hypersurface in $\Pj^n$. By
Proposition~\ref{prop:cones}, the quotient on the left above equals
\[
 \left(1-\frac1\delta\right)
 \frac{\pi_{n-2}(2)}{2^{n-2}}
 =\left(1-\frac1\delta\right)(2-2^{-(n-2)}).
\]
As odd $\delta$ tends to infinity these values tend to the asserted
constant, so no smaller one can work.

The supremum is not attained in any finite degree. A singular scheme has $\delta>1$, and the sharper
bound \eqref{eq:main} retains the factor
$(1-1/\delta)^{d-s}<1$ in the preceding upper estimate. If one uses
$\delta(\delta-1)^{d-s}q^s$ as the denominator, that factor
disappears. The upper bound is the same, and the binary cones over
$\F_2$ in Example~\ref{ex:quadratic-sharp} then give equality in
every degree $\delta\geqslant2$.
\end{proof}

\bibliographystyle{smfplain}
\bibliography{references}
\end{document}